\documentclass[
 superscriptaddress,
 amsmath,amssymb,
 aps,
 floatfix,
 twocolumn,prl
]{revtex4-2}

\usepackage{graphicx}
\usepackage{dcolumn}
\usepackage{bm}
\usepackage{xcolor}
\usepackage{hyperref}
\usepackage[capitalise]{cleveref}

\begin{document}

\preprint{APS/123-QED}

\title{Completely Integrable Replicator Dynamics Associated to Competitive Networks}

\author{Joshua Paik}
 \email{joshdpaik@gmail.com}
 \affiliation{%
 Department of Mathematics.
 The Pennsylvania State University, 
 University Park, PA 16802
}%
 
\author{Christopher Griffin}%
\email{griffinch@psu.edu}
\affiliation{
	Applied Research Laboratory,
    The Pennsylvania State University,
    University Park, PA 16802
    }%

\date{\today}

\begin{abstract} 
The replicator equations are a family of ordinary differential equations that arise in evolutionary game theory, and are closely related to Lotka-Volterra. We produce an infinite family of replicator equations which are Liouville-Arnold integrable. We show this by explicitly providing conserved quantities and a Poisson structure. As a corollary, we classify \textbf{all} tournament replicators up to dimension 6 and most of dimension 7. As an application, we show that Fig. 1 of ``A competitive network theory of species diversity" by Allesina and Levine (Proc. Natl. Acad. Sci., 2011), produces quasiperiodic dynamics.
\end{abstract}
\maketitle

The generalized Lotka-Volterra equations are, 
\begin{equation*}
\dot{x}_i = \epsilon_i x_i + x_i\left(\sum \limits_{j=1}^n W_{ij} x_j\right)
\end{equation*}
where $\epsilon_i$ are the ``intrinsic growth rates" and $W$ the interaction matrix. Closely related are the ``linear" replicator equations, 
\begin{equation}
    \dot{x}_i = x_i\left[(Wx)_i - x^TWx\right]
\end{equation}
These equations are used to theoretically explore population interactions in competitive ecosystems \cite{allesina2011competitive}. When $\epsilon_i = 0$ and $W$ is skew symmetric, Lotka-Volterra and the replicator equations are the same, and this is known to be a Hamiltonian system of ordinary differential equations (ODEs) for the Hamiltonian $H(x) = x_1 + ... + x_n$ and the quadratic Poisson Bracket \cref{eqn:Bracket} (see \cite{G21}). It is an interesting open question to characterize for which skew-symmetric $W$, the induced dynamics are integrable in the Liouville-Arnold sense. As an application, we know by KAM theory that any small enough Hamiltonian perturbation of the model also produces a large set of quasiperiodic solutions. (\cite{arnol2013mathematical}) What this means is that the qualitative behavior of the ODEs discussed in the remainder of the paper persists after a tiny change to the equations.

To this end, the Lotka-Volterra equation of the form \begin{equation} \label{eq:circulant}
    \dot{x}_i = x_i\left(\sum \limits_{j=1}^n x_{i+j} - x_{i-j}\right), \text{ for } i \in [1,..., 2n+1],
\end{equation} where indices are taken mod $2n+1$ has been shown to be integrable by Itoh \cite{itoh1987integrals} and Bogoyavlensky  \cite{bogoyavlensky1988five} \cite{bogoyavlensky1988integrable} and also by Veselov-Shabbat \cite{VS93} (equation 12). The latter is surprising as, while both Itoh and Bogoyavlensky were studying Lotka-Volterra, Veselov-Shabbat came up with a similar model as a discrete form of the Korteweg–De Vries (KdV) equations. Itoh did this by explicitly providing conserved quantities and showing that these conserved quantities Poisson commute. Both Bogoyavlensky and Veselov-Shabbat produced Lax Pairs. This connection was first observed in \cite{evripidou2017integrable}.

Observe that we may write \cref{eq:circulant} as \begin{equation}\label{eq:replicator}
    \dot{x}_i = x_i\left[\sum \limits_{j=1}^{2n+1} W_{ij}(G)x_j\right]
\end{equation}
where $W(G) = W$ is the tournament adjacency matrix of a circulant tournament graph of rock-paper-scissors type. By a circulant tournament $G$ of rock-paper-scissors type (hereafter referred to as a circulant tournament), we mean the tournament graph (up to isomorphism) such that for node every $i$, the incoming edges are of the form $\{(i-j) \to i \text{ (mod } 2n+1):j \in [1,...,n] \} $ and the outgoing edges are of the form $\{i \to (i+j) \text{ (mod } 2n+1):j \in [1,...,n] \}$. By tournament adjacency matrix, we mean the skew symmetric matrix $$W_{ij}(G) = \begin{cases}
    +1 & \text{if } j \to i \in G\\
    -1 & \text{if } i \to j \in G\\
    0 & \text{otherwise.}
\end{cases}$$

This is nice because we represent the ODE by a graph that encodes the equations. Hence, the \textit{combinatorics} of a given graph encodes the \textit{dynamics} of the corresponding differential equations. For example, the Volterra-Kac-van Moerbecke lattice is given by equation \cref{eq:replicator} when $G$ is a directed $n$-cycle. This was shown to be integrable by Kac-van Moerbecke \cite{kac1975explicitly} as well as Moser \cite{moser1975finitely}. 

\begin{figure}
    \centering
    \includegraphics[width=0.4\textwidth]{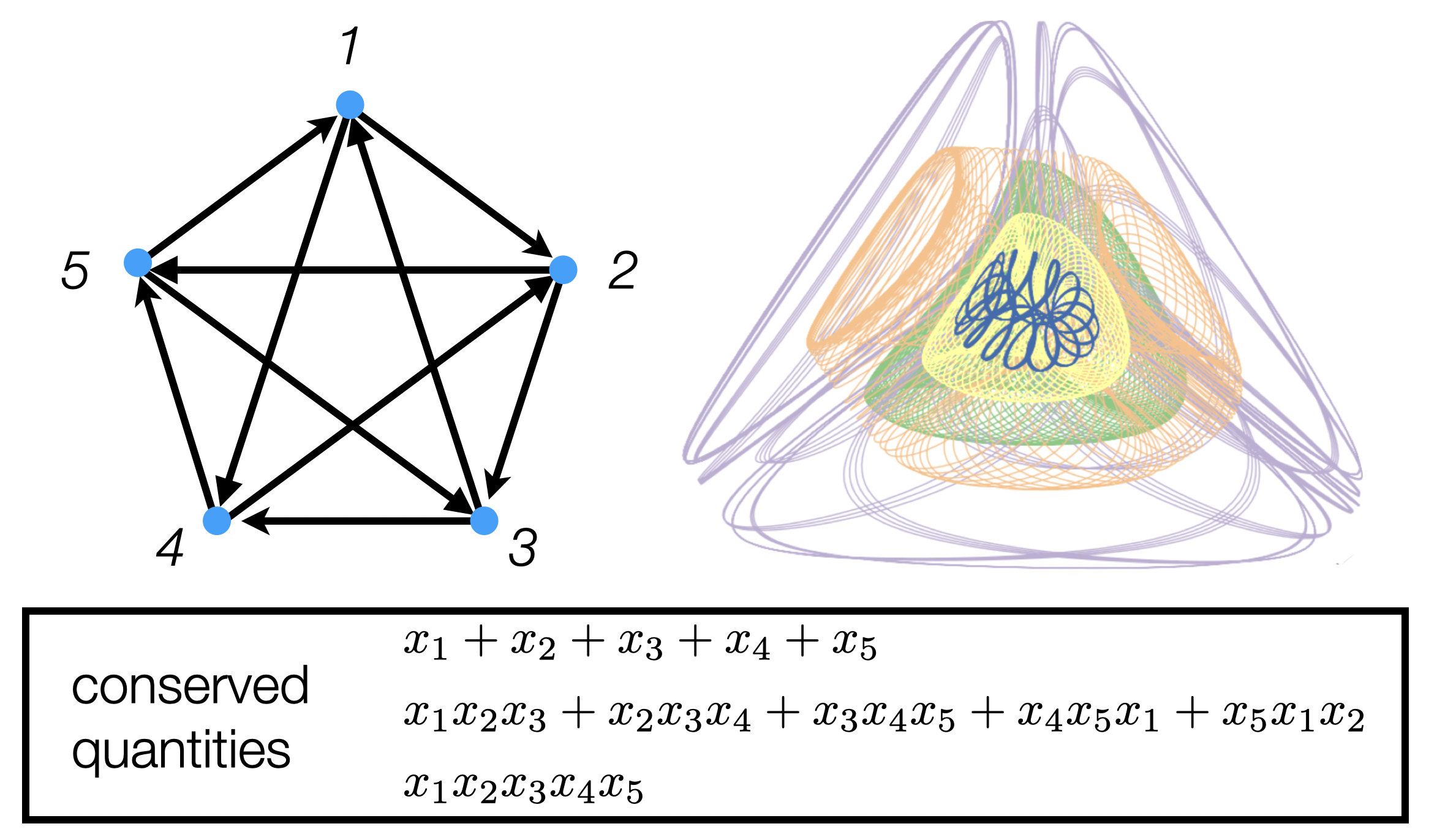}
    \caption{A circulant tournament of order 5 is shown with labelled nodes. Its dynamics is shown to the right, and one sees a quasiperiodic foliation by tori. The conserved quantities of Itoh are shown.}
    \label{fig:1}
\end{figure}

Our first result is to repackage Itoh's conserved quantity result with the philosophy that combinatorics informs dynamics. In short, we notice that the polynomials defining the conserved quantities correspond to subgraphs of the circulant tournament $G.$ Let $C_k(G)$ be the set of subsets of $\{1,...,2n+1\}$ such that $\sigma \in C_k(G)$ induces a circulant tournament of order $k$. Then, for every odd $k$, \begin{equation}
    F_k(G)(x) = \sum \limits_{\sigma \in C_k(G)} \prod \limits_{j \in \sigma} x_j
    \label{eq:conservations}
\end{equation}
is a conserved quantity of \cref{eq:replicator}. The proof is the same as Itoh's, but instead of a careful accounting of indices as Itoh, one does a careful accounting of paths. 

The quadratic bracket \begin{equation}
    \{F,H\}_W = 
\sum_{i < j} W_{ij}x_ix_j\left(\frac{\partial F}{\partial x_i} \frac{\partial H}{\partial x_j} - \frac{\partial F}{\partial x_j} \frac{\partial H}{\partial x_i} \right).
\label{eqn:Bracket}
\end{equation}
is a Poisson bracket such that the conserved quantities of \cref{eq:conservations} Poisson commute. This was shown by Itoh as well in his stunning paper \cite{itoh2008combinatorial}. Under this bracket, $F_{2n+1} = x_1x_2...x_{2n+1}$ is the Casimir. By the Liouville-Arnold theorem, the dynamics of \cref{eq:replicator} for circulant tournament $G$ are quasiperiodic and foliate the phase space by tori. It is interesting to note this bracket is identical to the one used by Ovsienko, Schwarz, and Tabachnikov (albeit for a different $W$) in their proof of the integrability of the Pentagram map \cite{ovsienko2013liouville}.

The upshot of this graphical interpretation of Itoh's result, is that we can produced conserved quantities for \cref{eq:replicator} for a larger class of graphs. We call this the embedding procedure. This procedure is heavily inspired by the ``cyclic union" technique developed by Curto and collaborators \cite{parmelee2022core,parmelee2022sequential,morrison2019predicting,santander2022nerve}, in the context of neural network dynamics, but adapted to the context of circulant tournaments. 
\begin{figure}[t]
    \centering
    \includegraphics[width=0.4\textwidth]{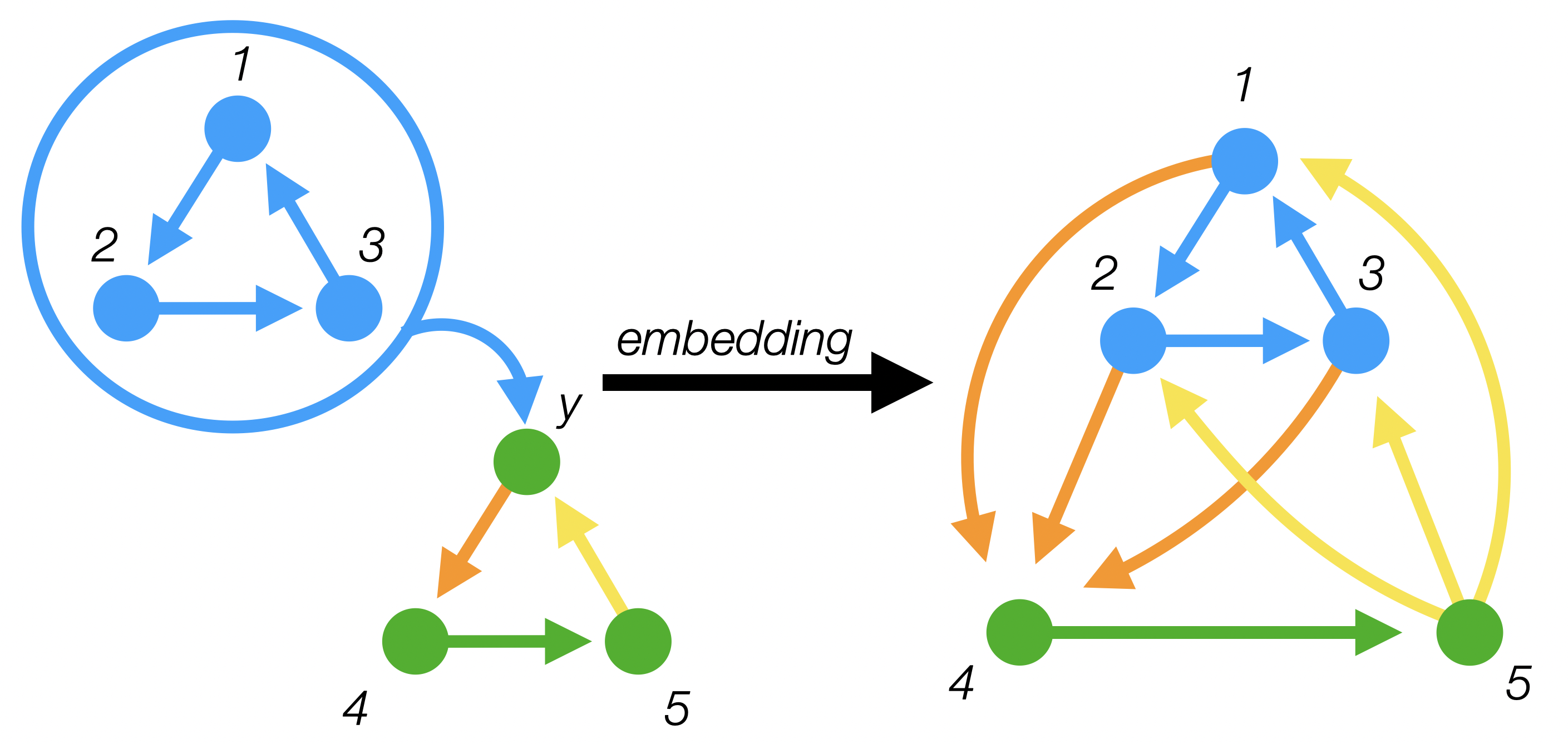}

    \includegraphics[width = 0.4\textwidth]{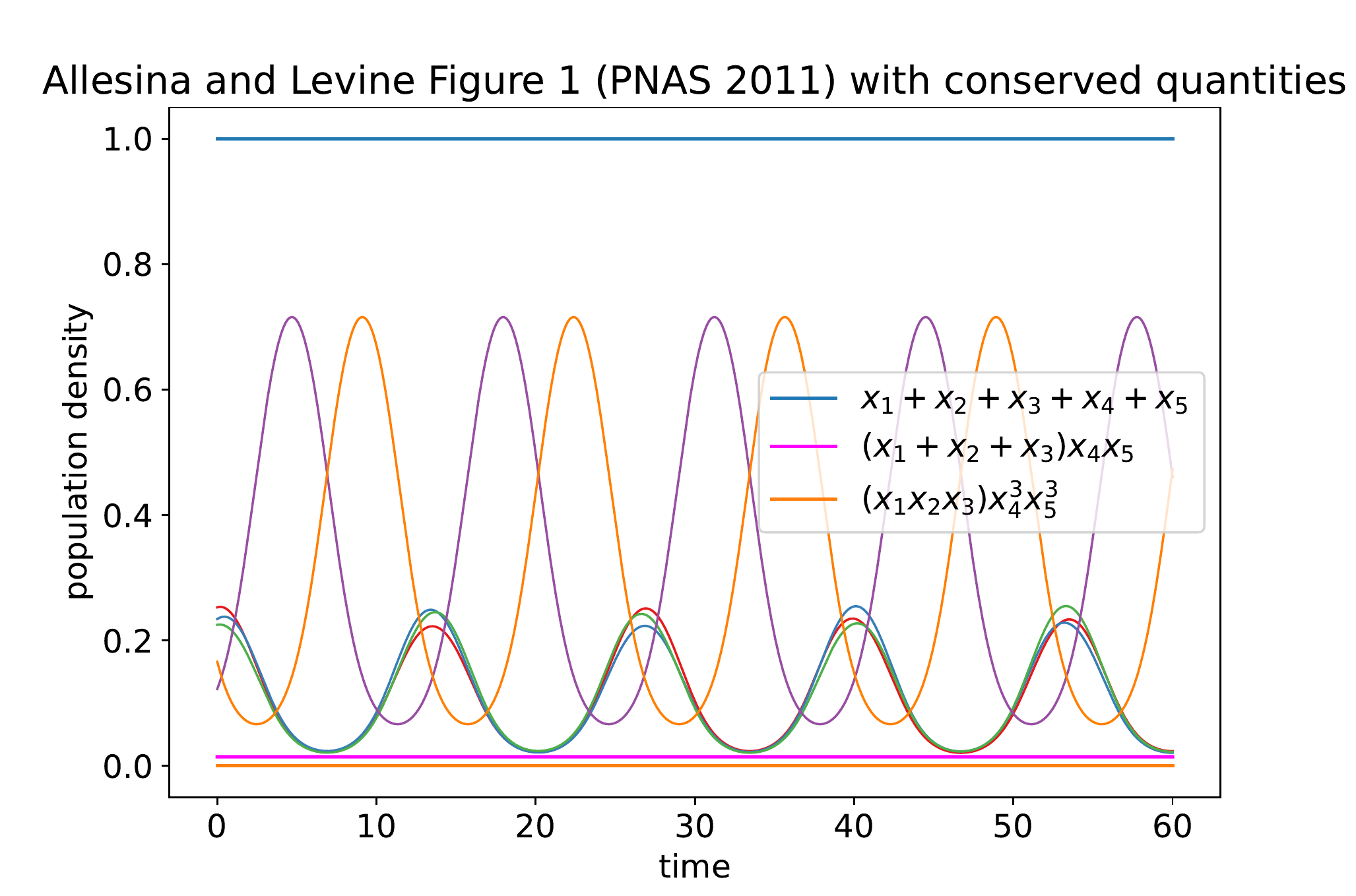}
    \caption{We depict how to embed a 3-cycle into a single vertex of another 3-cycle. In the bottom figure, we depict solutions and conserved quantities for the corresponding replicator dynamics. Interestingly, this figure without the conserved quantities were already in the literature, Figure 1 of \cite{allesina2011competitive}. By our conserved quantities result and Liouville-Arnold theorem, we know that the dynamics are quasiperiodic. }
    \label{fig:allesina}
\end{figure}

To be more specific about embedding, let $G$ be a circulant tournament on $2n+1$ nodes. Let $H = \left(H^1,...,H^{2n+1}\right)$ be an ordered list of graphs such that the number of vertices $H^i$ is $2n_i+1$ for some $n_i \in \{0,1,2,\dots\}$. Then we embed $H$ into $G$, denoted $G \hookleftarrow H$, by replacing every node $i \in G$ with $H^i$, and forming an edge $h^i_q \to h^i_v$ if $q \to v$ was an edge in $G$. (See \cref{fig:allesina}.) We can repeat this to produce a chain of embeddings, like Russian Matryoshka dolls. When $G$ is a 3 cycle and $H$ is a 3-cycle and two singleton graphs (single vertices), then $G \hookleftarrow H$ is produced in \cref{fig:allesina}. Interestingly, this graph is precisely Fig. 1 of \cite{allesina2011competitive}. The conserved quantities are given in the figure, and as a corollary, Fig. 1 of \cite{allesina2011competitive} is Liouville-Arnold integrable and quasiperiodic. 

The power to produce conserved quantity results for embedded circulant tournaments comes from a simple observation. Suppose we embed a single graph $H$ into a node $y$ of $G$, which we denote $J = G \hookleftarrow H$. The degree $k$ conserved quantity $h_k(x)(t)$ for the replicator dynamics induced by $H$, is \textit{not} a conserved quantity for the replicator dynamics of $J$, when we appropriately substitute variables of $J$ into $h_k$. However, $h_k(x)(t)$ in fact equals $y(t)^k$ where $y$ is seen as a solution of the dynamics for the corresponding graph $G$. This can be shown by a direct computation. Therefore, if $f(y, x_2,...,x_{2n+1})$ was a conserved quantity for $G$ alone, then $f(h_k, x_2^k,...,x_{2n+1}^k)$ is a conserved quantity for the embedded graph. Also, note that if $f_k$ is a conserved quantity, then $\alpha f_k^n$ is a conserved quantity for any constant $\alpha$ and integer $n$. 
\begin{figure*}
\centering
\includegraphics[width=0.8\textwidth]{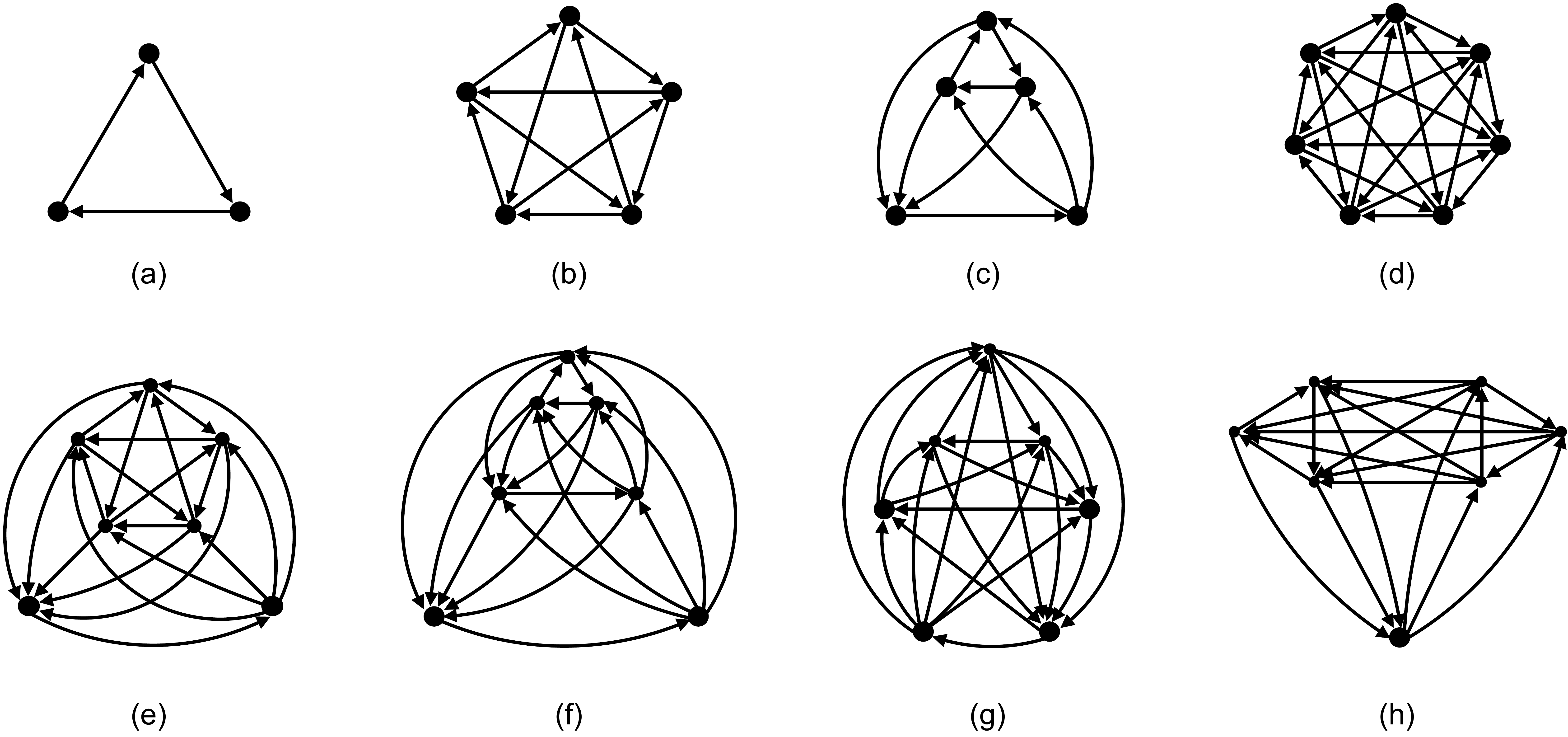}
\caption{Depicted are all tournaments up to dimension 7 which are produced by starting with 3, 5, and 7 circulant (rock-paper-scissor-like) tournaments, and creating appropriate embeddings. The conserved quantities are conveniently listed in the appendix. All of these have completely integrable dynamics. As a corollary of this result, we can classify all zero-sum classical tournament replicators up to dimension 5 and all but seven in dimension 7. The 7 graphs which are not classified appear, numerically, to have chaotic (positive Lyapunov exponent) dynamics.}
\label{fig:IntegrableGraphs}
\end{figure*}

We now proceed to produce conserved quantities for embedded graphs. In the simplest case, we produce conserved quantities of $J = G \hookleftarrow H$ where $G$ is a circulant tournament of order $2n+1$ and $H$ is a list containing precisely one circulant tournament of order $2m+1$, $m > 0$ and the rest, singletons. Without loss of generality, suppose the  circulant tournament to be embedded is $H^1$ with vertices $\{h^1_i\}_i$. Then $J$ is a tournament of order $2(n+m) + 1$. We need precisely $n+m + 1$ conserved quantities in order to invoke the Liouville-Arnold theorem. If $g_k(x_1,...,x_{2n+1})$ is a conserved quantity of $G$, then $g_k(\sum_{h_i^1\in H^1} h_i^1, x_2,...,x_{2n+1})$ is a conserved quantity of $J$ for every odd $k \in [1,...,n]$. This gives $2n+1$ conserved quantities, which we call the \textit{outer symmetries}. Let $g_{2m+1} = x_1\prod_{i=2}^n x_i^{x_i^*}$ Now, let $h_k$ be the degree $k$ conserved quantity of $H$. Then the second conjugacy result implies $h_k\left(\prod_{i=2}^n x_i^{x_i^*}\right)^k$ is a conserved quantity. This gives $m+1$ conserved quantities. We call these the \textit{inner symmetries}. But notice, one is repeated, so this process produces $n+m+1$ conserved quantities of $J$.

As an example, we compute the conserved quantities for the graph in \cref{fig:allesina}. The outer symmetries of $J = G \hookleftarrow H$ are, 
\begin{align*}
    C_1 &= (x_1 + x_2 + x_3) + x_4 + x_5 \\
    C_3 &= (x_1+x_2+x_3)x_4x_5,
\end{align*}
which are produced by substituting $x_1 + x_2 + x_3$ for $y$. 

The inner symmetries of $J$ are, 
\begin{align*}
    C_3 &= (x_1+x_2+x_3)x_4x_5\\
    C_9 &= (x_1x_2x_3)x_4^3x_5^3.
\end{align*}
Therefore, we have three conserved quantities -- $C_1, C_3$ and $C_9$. It follows from a small variation on Itoh's result that these conserved quantities Poisson commute. Therefore, we produced the conserved quantities numerically depicted in \cref{fig:allesina}. 

When we embed more than one non-trivial circulant tournament into G, we follow the same procedure, except we sum over all possible conserved quantities produced from the same degree. (There is a technicality -- we should sum over at most $|G|$ polynomials of the same order.) This is a result of the (almost trivial) observation that constant multiples of a conserved quantity are still conserved. For example, in \cref{fig:IntegrableGraphs}(h), label the vertices of one of the three-cycles (i.e., rock-paper-scissors tournaments) as $1,2,3$ embedded in $y_1$ and the vertices of the second one as $4,5,6$ embedded in $y_2$. Then the degree 7 conserved quantity is,
\begin{equation*}
(x_1+x_2+x_3)^3(x_4x_5x_6)x_7 + (x_1x_2x_3)(x_4+x_5+x_6)^3x_7.
\end{equation*}
However, notice that both terms are equivalent to $y_1^3y_2^3x_7^3$, when the $y_i$ are from the original graph, so their sum can be mapped to $2y_1^3y_2^3x_7^3$

This embedding procedure and two folk theorems are enough to classify the dynamics of all zero-sum tournament based replicators up to dimension 6 and all but 7 in dimension 7. The remaining 7 graphs in dimension 7 with a full support fixed point appear to have chaotic (positive Lyapunov exponent) dynamics, and elude classification. The two folk theorems are as follows. First, the size of the support of a fixed point must be odd \cite{mccarthy1996determinants}. Second, the time average of an orbit, $\lim \limits_{T \to \infty} \frac{1}{T} \int \limits_0^T x(t) dt$ equals one of the fixed points. As the dynamics are positive, these two results combined mean we must only analyze tournaments which admit a full support fixed point -- which only happens in odd dimensions. If a graph does not admit a full support fixed point, then it necessarily reduces to a smaller, odd dimensional system.

We now start the classification. We use Brendan McKay's database of tournaments \cite{mckay}. In dimension 3, there is precisely one tournament with a full support fixed point. This is the three cycle. In dimension 5, there are 12 tournaments. Of those, precisely two have full support fixed points, they are shown in \cref{fig:IntegrableGraphs} b and c. The graph in \cref{fig:IntegrableGraphs}(b) is the 5 RPS, and the graph in \cref{fig:IntegrableGraphs}(c) is the graph shown in \cref{fig:allesina} -- both are completely integrable. The remaining 10 tournament graphs in dimension 5 either reduce to a 3-cycle or a singleton fixed point. In dimension 7, there are 456 tournaments. Of those, 12 have a full support fixed point. Of those with a full support fixed point, we can prove that 5 are completely integrable -- they are shown in \cref{fig:IntegrableGraphs} d - h. 

\begin{figure}[htbp]
\centering
\includegraphics[width=0.4\columnwidth]{graph10.pdf}\quad \includegraphics[width=0.5\columnwidth]{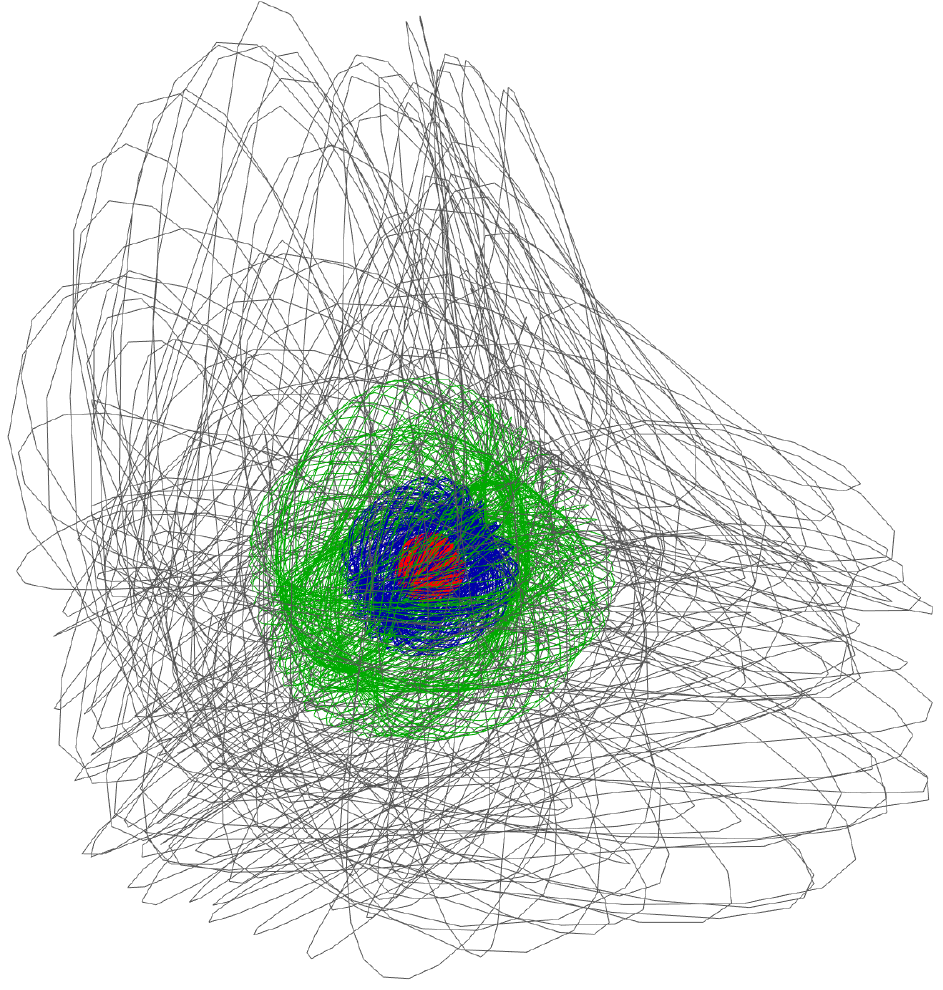}
\caption{The graph on the left produces replicator dynamics with phase portrait on the right. One observes quasiperiodic islands near the elliptic fixed point, and chaotic orbits the further out one goes. Numerical experiments show positivity of Lyapunov Exponents.}
\label{fig:Graph10}
\end{figure}

Of the 7 graphs remaining in dimension 7 that were not classified, only the graph in \cref{fig:Graph10} has a non-trivial graph automorphism group. The dynamics are volume preserving. It is interesting because, near the elliptic full support fixed point, we see quasiperiodic dynamics, however, as we move further from the fixed point, the dynamics become more and more chaotic. This is shown in the phase portraits, with elliptical orbits colored, and chaotic orbits in black. This agrees with numerical experiments on the Lyapunov spectra, which are symmetric around 0. This is a typical picture in symplectic, Hamiltonian, volume preserving dynamics, elliptical islands floating in ergodic seas.

It seems worthwhile to comment on what happens when $W$ has the same arrangement of positive and negative entries as a tournament, but we weaken the condition that entries are $\pm 1$. Here is an example which suggests some tractability. If instead of the typical adjacency of the graph in \cref{fig:allesina}, we had $$W = \begin{bmatrix}
     0 &-5 &+1 &-1 &+1\\
    +5 & 0 &-1 &-1 &+1\\
    -1 &+1 &0 &-1 &+1\\
    +1 &+1 &+1 & 0 &-1\\
    -1 &-1 &-1 &+1& 0
\end{bmatrix}$$ 
then the conserved quantities for the corresponding replicator dynamics are \begin{align*}
f_1 &= x_1+x_2+x_3+x_4+x_5\\
f_3& = (x_1+x_2+x_3)x_4x_5 \\
f_{21}& = x_1x_2x_3^5x_4^7x_5^7
\end{align*}
However, one can easily find skew symmetric $W$ that appear numerically to have elliptic islands in a sea of ergodicity, even in dimension 5. How to deal with this is unclear and represents an area of future research. Likewise, merging the work presented in \cite{E19,charalambides2015generalized,C19, evripidou2017integrable} with the embedding approach presented in this letter may also lead to novel results in integrability. However, care must be taken in these situations because chaos seems to easily arise from small perturbations of these systems.

Related work on integrable systems arising from Lotka-Volterra (LV) dynamics can is discussed by Charalambides et al. \cite{charalambides2015generalized}. While considering a broad class of dynamics, \cite{charalambides2015generalized} does not consider dynamics arising from embedded tournament structures. In particular, the examples provided lack tournament or embedding structure. Moreover, their approach uses Lax pair arguments, while the results in this paper focus on LV dynamics arising from the replicator equation applied to tournaments constructed by graph embeddings. Our approach is fundamentally combinatorial. Work by Evripidou, Kassotakis and Vanhaecke \cite{evripidou2017integrable}, provides an alternative proof of the results in \cite{VS93, itoh1987integrals} using a deformed Casimir, but again uses a complete tournament structure, rather than a tournament arising from a graph embedding. Additional work by these authors in \cite{E19} is fascinating and appears to extend the integrability results of the Volterra lattice to graphs with additional cyclic symmetries. Unlike the results presented, they do not consider graphs (tournaments) that result from graph embeddings. However, a combination of the work in \cite{E19} and the work presented in this paper is a potential area of future research. Finally, \cite{C19} studies the problem presented in \cite{VS93,itoh1987integrals} but removes the criteria on the constant terms that are added to the dynamics. These dynamics do not arise from graph embeddings, as considered here.

This work was supported by NSF CMMI-1932991. We thank Sergei Tabachnikov and Sasha Veselov for pointing out the connection to the Dressing Chain. We thank the reviewers for helpful comments which significantly improved this paper.

\bibliography{bibliography}

\begin{thebibliography}{21}%
\makeatletter
\providecommand \@ifxundefined [1]{%
 \@ifx{#1\undefined}
}%
\providecommand \@ifnum [1]{%
 \ifnum #1\expandafter \@firstoftwo
 \else \expandafter \@secondoftwo
 \fi
}%
\providecommand \@ifx [1]{%
 \ifx #1\expandafter \@firstoftwo
 \else \expandafter \@secondoftwo
 \fi
}%
\providecommand \natexlab [1]{#1}%
\providecommand \enquote  [1]{``#1''}%
\providecommand \bibnamefont  [1]{#1}%
\providecommand \bibfnamefont [1]{#1}%
\providecommand \citenamefont [1]{#1}%
\providecommand \href@noop [0]{\@secondoftwo}%
\providecommand \href [0]{\begingroup \@sanitize@url \@href}%
\providecommand \@href[1]{\@@startlink{#1}\@@href}%
\providecommand \@@href[1]{\endgroup#1\@@endlink}%
\providecommand \@sanitize@url [0]{\catcode `\\12\catcode `\$12\catcode
  `\&12\catcode `\#12\catcode `\^12\catcode `\_12\catcode `\%12\relax}%
\providecommand \@@startlink[1]{}%
\providecommand \@@endlink[0]{}%
\providecommand \url  [0]{\begingroup\@sanitize@url \@url }%
\providecommand \@url [1]{\endgroup\@href {#1}{\urlprefix }}%
\providecommand \urlprefix  [0]{URL }%
\providecommand \Eprint [0]{\href }%
\providecommand \doibase [0]{https://doi.org/}%
\providecommand \selectlanguage [0]{\@gobble}%
\providecommand \bibinfo  [0]{\@secondoftwo}%
\providecommand \bibfield  [0]{\@secondoftwo}%
\providecommand \translation [1]{[#1]}%
\providecommand \BibitemOpen [0]{}%
\providecommand \bibitemStop [0]{}%
\providecommand \bibitemNoStop [0]{.\EOS\space}%
\providecommand \EOS [0]{\spacefactor3000\relax}%
\providecommand \BibitemShut  [1]{\csname bibitem#1\endcsname}%
\let\auto@bib@innerbib\@empty
\bibitem [{\citenamefont {Allesina}\ and\ \citenamefont
  {Levine}(2011)}]{allesina2011competitive}%
  \BibitemOpen
  \bibfield  {author} {\bibinfo {author} {\bibfnamefont {S.}~\bibnamefont
  {Allesina}}\ and\ \bibinfo {author} {\bibfnamefont {J.~M.}\ \bibnamefont
  {Levine}},\ }\bibfield  {title} {\bibinfo {title} {A competitive network
  theory of species diversity},\ }\href@noop {} {\bibfield  {journal} {\bibinfo
   {journal} {Proceedings of the National Academy of Sciences}\ }\textbf
  {\bibinfo {volume} {108}},\ \bibinfo {pages} {5638} (\bibinfo {year}
  {2011})}\BibitemShut {NoStop}%
\bibitem [{\citenamefont {Griffin}(2021)}]{G21}%
  \BibitemOpen
  \bibfield  {author} {\bibinfo {author} {\bibfnamefont {C.}~\bibnamefont
  {Griffin}},\ }\bibfield  {title} {\bibinfo {title} {The replicator dynamics
  of zero-sum games arise from a novel {Poisson} algebra},\ }\href@noop {}
  {\bibfield  {journal} {\bibinfo  {journal} {Chaos, Solitons \& Fractals}\
  }\textbf {\bibinfo {volume} {153}},\ \bibinfo {pages} {111508} (\bibinfo
  {year} {2021})}\BibitemShut {NoStop}%
\bibitem [{\citenamefont {Arnol'd}(2013)}]{arnol2013mathematical}%
  \BibitemOpen
  \bibfield  {author} {\bibinfo {author} {\bibfnamefont {V.~I.}\ \bibnamefont
  {Arnol'd}},\ }\href@noop {} {\emph {\bibinfo {title} {Mathematical methods of
  classical mechanics}}},\ Vol.~\bibinfo {volume} {60}\ (\bibinfo  {publisher}
  {Springer Science \& Business Media},\ \bibinfo {year} {2013})\BibitemShut
  {NoStop}%
\bibitem [{\citenamefont {Itoh}(1987)}]{itoh1987integrals}%
  \BibitemOpen
  \bibfield  {author} {\bibinfo {author} {\bibfnamefont {Y.}~\bibnamefont
  {Itoh}},\ }\bibfield  {title} {\bibinfo {title} {Integrals of a
  lotka-volterra system of odd number of variables},\ }\href@noop {} {\bibfield
   {journal} {\bibinfo  {journal} {Progress of theoretical physics}\ }\textbf
  {\bibinfo {volume} {78}},\ \bibinfo {pages} {507} (\bibinfo {year}
  {1987})}\BibitemShut {NoStop}%
\bibitem [{\citenamefont
  {Bogoyavlensky}(1988{\natexlab{a}})}]{bogoyavlensky1988five}%
  \BibitemOpen
  \bibfield  {author} {\bibinfo {author} {\bibfnamefont {O.}~\bibnamefont
  {Bogoyavlensky}},\ }\bibfield  {title} {\bibinfo {title} {Five constructions
  of integrable dynamical systems connected with the korteweg-de vries
  equation},\ }\href@noop {} {\bibfield  {journal} {\bibinfo  {journal} {Acta
  Applicandae Mathematica}\ }\textbf {\bibinfo {volume} {13}},\ \bibinfo
  {pages} {227} (\bibinfo {year} {1988}{\natexlab{a}})}\BibitemShut {NoStop}%
\bibitem [{\citenamefont
  {Bogoyavlensky}(1988{\natexlab{b}})}]{bogoyavlensky1988integrable}%
  \BibitemOpen
  \bibfield  {author} {\bibinfo {author} {\bibfnamefont {O.}~\bibnamefont
  {Bogoyavlensky}},\ }\bibfield  {title} {\bibinfo {title} {Integrable
  discretizations of the kdv equation},\ }\href@noop {} {\bibfield  {journal}
  {\bibinfo  {journal} {Physics Letters A}\ }\textbf {\bibinfo {volume}
  {134}},\ \bibinfo {pages} {34} (\bibinfo {year}
  {1988}{\natexlab{b}})}\BibitemShut {NoStop}%
\bibitem [{\citenamefont {Veselov}\ and\ \citenamefont {Shabat}(1993)}]{VS93}%
  \BibitemOpen
  \bibfield  {author} {\bibinfo {author} {\bibfnamefont {A.~P.}\ \bibnamefont
  {Veselov}}\ and\ \bibinfo {author} {\bibfnamefont {A.~B.}\ \bibnamefont
  {Shabat}},\ }\bibfield  {title} {\bibinfo {title} {Dressing chains and
  spectral theory of the schrodinger operator},\ }\href@noop {} {\bibfield
  {journal} {\bibinfo  {journal} {Funktsional'nyi Analiz i ego Prilozheniya}\
  }\textbf {\bibinfo {volume} {27}},\ \bibinfo {pages} {1} (\bibinfo {year}
  {1993})}\BibitemShut {NoStop}%
\bibitem [{\citenamefont {Evripidou}\ \emph {et~al.}(2017)\citenamefont
  {Evripidou}, \citenamefont {Kassotakis},\ and\ \citenamefont
  {Vanhaecke}}]{evripidou2017integrable}%
  \BibitemOpen
  \bibfield  {author} {\bibinfo {author} {\bibfnamefont {C.~A.}\ \bibnamefont
  {Evripidou}}, \bibinfo {author} {\bibfnamefont {P.}~\bibnamefont
  {Kassotakis}},\ and\ \bibinfo {author} {\bibfnamefont {P.}~\bibnamefont
  {Vanhaecke}},\ }\bibfield  {title} {\bibinfo {title} {Integrable deformations
  of the bogoyavlenskij--itoh lotka--volterra systems},\ }\href@noop {}
  {\bibfield  {journal} {\bibinfo  {journal} {Regular and Chaotic Dynamics}\
  }\textbf {\bibinfo {volume} {22}},\ \bibinfo {pages} {721} (\bibinfo {year}
  {2017})}\BibitemShut {NoStop}%
\bibitem [{\citenamefont {Kac}\ and\ \citenamefont {van
  Moerbeke}(1975)}]{kac1975explicitly}%
  \BibitemOpen
  \bibfield  {author} {\bibinfo {author} {\bibfnamefont {M.}~\bibnamefont
  {Kac}}\ and\ \bibinfo {author} {\bibfnamefont {P.}~\bibnamefont {van
  Moerbeke}},\ }\bibfield  {title} {\bibinfo {title} {On an explicitly soluble
  system of nonlinear differential equations related to certain toda
  lattices},\ }\href@noop {} {\bibfield  {journal} {\bibinfo  {journal}
  {Advances in Mathematics}\ }\textbf {\bibinfo {volume} {16}},\ \bibinfo
  {pages} {160} (\bibinfo {year} {1975})}\BibitemShut {NoStop}%
\bibitem [{\citenamefont {Moser}(1975)}]{moser1975finitely}%
  \BibitemOpen
  \bibfield  {author} {\bibinfo {author} {\bibfnamefont {J.}~\bibnamefont
  {Moser}},\ }\bibfield  {title} {\bibinfo {title} {Finitely many mass points
  on the line under the influence of an exponential potential--an integrable
  system},\ }\href@noop {} {\bibfield  {journal} {\bibinfo  {journal}
  {Dynamical systems, theory and applications}\ ,\ \bibinfo {pages} {467}}
  (\bibinfo {year} {1975})}\BibitemShut {NoStop}%
\bibitem [{\citenamefont {Itoh}(2008)}]{itoh2008combinatorial}%
  \BibitemOpen
  \bibfield  {author} {\bibinfo {author} {\bibfnamefont {Y.}~\bibnamefont
  {Itoh}},\ }\bibfield  {title} {\bibinfo {title} {A combinatorial method for
  the vanishing of the poisson brackets of an integrable lotka--volterra
  system},\ }\href@noop {} {\bibfield  {journal} {\bibinfo  {journal} {Journal
  of Physics A: Mathematical and Theoretical}\ }\textbf {\bibinfo {volume}
  {42}},\ \bibinfo {pages} {025201} (\bibinfo {year} {2008})}\BibitemShut
  {NoStop}%
\bibitem [{\citenamefont {Ovsienko}\ \emph {et~al.}(2013)\citenamefont
  {Ovsienko}, \citenamefont {Schwartz},\ and\ \citenamefont
  {Tabachnikov}}]{ovsienko2013liouville}%
  \BibitemOpen
  \bibfield  {author} {\bibinfo {author} {\bibfnamefont {V.}~\bibnamefont
  {Ovsienko}}, \bibinfo {author} {\bibfnamefont {R.~E.}\ \bibnamefont
  {Schwartz}},\ and\ \bibinfo {author} {\bibfnamefont {S.}~\bibnamefont
  {Tabachnikov}},\ }\bibfield  {title} {\bibinfo {title} {Liouville--arnold
  integrability of the pentagram map on closed polygons},\ }\href@noop {}
  {\bibfield  {journal} {\bibinfo  {journal} {Duke Mathematical Journal}\
  }\textbf {\bibinfo {volume} {162}},\ \bibinfo {pages} {2149} (\bibinfo {year}
  {2013})}\BibitemShut {NoStop}%
\bibitem [{\citenamefont {Parmelee}\ \emph
  {et~al.}(2022{\natexlab{a}})\citenamefont {Parmelee}, \citenamefont {Moore},
  \citenamefont {Morrison},\ and\ \citenamefont {Curto}}]{parmelee2022core}%
  \BibitemOpen
  \bibfield  {author} {\bibinfo {author} {\bibfnamefont {C.}~\bibnamefont
  {Parmelee}}, \bibinfo {author} {\bibfnamefont {S.}~\bibnamefont {Moore}},
  \bibinfo {author} {\bibfnamefont {K.}~\bibnamefont {Morrison}},\ and\
  \bibinfo {author} {\bibfnamefont {C.}~\bibnamefont {Curto}},\ }\bibfield
  {title} {\bibinfo {title} {Core motifs predict dynamic attractors in
  combinatorial threshold-linear networks},\ }\href@noop {} {\bibfield
  {journal} {\bibinfo  {journal} {PloS one}\ }\textbf {\bibinfo {volume}
  {17}},\ \bibinfo {pages} {e0264456} (\bibinfo {year}
  {2022}{\natexlab{a}})}\BibitemShut {NoStop}%
\bibitem [{\citenamefont {Parmelee}\ \emph
  {et~al.}(2022{\natexlab{b}})\citenamefont {Parmelee}, \citenamefont
  {Alvarez}, \citenamefont {Curto},\ and\ \citenamefont
  {Morrison}}]{parmelee2022sequential}%
  \BibitemOpen
  \bibfield  {author} {\bibinfo {author} {\bibfnamefont {C.}~\bibnamefont
  {Parmelee}}, \bibinfo {author} {\bibfnamefont {J.~L.}\ \bibnamefont
  {Alvarez}}, \bibinfo {author} {\bibfnamefont {C.}~\bibnamefont {Curto}},\
  and\ \bibinfo {author} {\bibfnamefont {K.}~\bibnamefont {Morrison}},\
  }\bibfield  {title} {\bibinfo {title} {Sequential attractors in combinatorial
  threshold-linear networks},\ }\href@noop {} {\bibfield  {journal} {\bibinfo
  {journal} {SIAM Journal on Applied Dynamical Systems}\ }\textbf {\bibinfo
  {volume} {21}},\ \bibinfo {pages} {1597} (\bibinfo {year}
  {2022}{\natexlab{b}})}\BibitemShut {NoStop}%
\bibitem [{\citenamefont {Morrison}\ and\ \citenamefont
  {Curto}(2019)}]{morrison2019predicting}%
  \BibitemOpen
  \bibfield  {author} {\bibinfo {author} {\bibfnamefont {K.}~\bibnamefont
  {Morrison}}\ and\ \bibinfo {author} {\bibfnamefont {C.}~\bibnamefont
  {Curto}},\ }\bibfield  {title} {\bibinfo {title} {Predicting neural network
  dynamics via graphical analysis},\ }in\ \href@noop {} {\emph {\bibinfo
  {booktitle} {Algebraic and Combinatorial Computational Biology}}}\ (\bibinfo
  {publisher} {Elsevier},\ \bibinfo {year} {2019})\ pp.\ \bibinfo {pages}
  {241--277}\BibitemShut {NoStop}%
\bibitem [{\citenamefont {Santander}\ \emph {et~al.}(2022)\citenamefont
  {Santander}, \citenamefont {Ebli}, \citenamefont {Patania}, \citenamefont
  {Sanderson}, \citenamefont {Burtscher}, \citenamefont {Morrison},\ and\
  \citenamefont {Curto}}]{santander2022nerve}%
  \BibitemOpen
  \bibfield  {author} {\bibinfo {author} {\bibfnamefont {D.~E.}\ \bibnamefont
  {Santander}}, \bibinfo {author} {\bibfnamefont {S.}~\bibnamefont {Ebli}},
  \bibinfo {author} {\bibfnamefont {A.}~\bibnamefont {Patania}}, \bibinfo
  {author} {\bibfnamefont {N.}~\bibnamefont {Sanderson}}, \bibinfo {author}
  {\bibfnamefont {F.}~\bibnamefont {Burtscher}}, \bibinfo {author}
  {\bibfnamefont {K.}~\bibnamefont {Morrison}},\ and\ \bibinfo {author}
  {\bibfnamefont {C.}~\bibnamefont {Curto}},\ }\bibfield  {title} {\bibinfo
  {title} {Nerve theorems for fixed points of neural networks},\ }in\
  \href@noop {} {\emph {\bibinfo {booktitle} {Research in Computational
  Topology 2}}}\ (\bibinfo  {publisher} {Springer},\ \bibinfo {year} {2022})\
  pp.\ \bibinfo {pages} {129--165}\BibitemShut {NoStop}%
\bibitem [{\citenamefont {McCarthy}\ and\ \citenamefont
  {Benjamin}(1996)}]{mccarthy1996determinants}%
  \BibitemOpen
  \bibfield  {author} {\bibinfo {author} {\bibfnamefont {C.~A.}\ \bibnamefont
  {McCarthy}}\ and\ \bibinfo {author} {\bibfnamefont {A.~T.}\ \bibnamefont
  {Benjamin}},\ }\bibfield  {title} {\bibinfo {title} {Determinants of the
  tournaments},\ }\href@noop {} {\bibfield  {journal} {\bibinfo  {journal}
  {Mathematics Magazine}\ }\textbf {\bibinfo {volume} {69}},\ \bibinfo {pages}
  {133} (\bibinfo {year} {1996})}\BibitemShut {NoStop}%
\bibitem [{\citenamefont {McKay}()}]{mckay}%
  \BibitemOpen
  \bibfield  {author} {\bibinfo {author} {\bibfnamefont {B.}~\bibnamefont
  {McKay}},\ }\href {https://users.cecs.anu.edu.au/~bdm/data/digraphs.html}
  {\bibinfo {title} {Digraphs}},\ \bibinfo {note} {accessed:
  2022-10-27}\BibitemShut {NoStop}%
\bibitem [{E19(2019)}]{E19}%
  \BibitemOpen
  \bibfield  {title} {\bibinfo {title} {Integrable reductions of the dressing
  chain},\ }\href {https://doi.org/10.3934/jcd.2019014} {\bibfield  {journal}
  {\bibinfo  {journal} {Journal of Computational Dynamics}\ }\textbf {\bibinfo
  {volume} {6}},\ \bibinfo {pages} {277} (\bibinfo {year} {2019})}\BibitemShut
  {NoStop}%
\bibitem [{\citenamefont {Charalambides}\ \emph {et~al.}(2015)\citenamefont
  {Charalambides}, \citenamefont {Damianou},\ and\ \citenamefont
  {Evripidou}}]{charalambides2015generalized}%
  \BibitemOpen
  \bibfield  {author} {\bibinfo {author} {\bibfnamefont {S.~A.}\ \bibnamefont
  {Charalambides}}, \bibinfo {author} {\bibfnamefont {P.~A.}\ \bibnamefont
  {Damianou}},\ and\ \bibinfo {author} {\bibfnamefont {C.~A.}\ \bibnamefont
  {Evripidou}},\ }\bibfield  {title} {\bibinfo {title} {On generalized volterra
  systems},\ }\href@noop {} {\bibfield  {journal} {\bibinfo  {journal} {Journal
  of Geometry and Physics}\ }\textbf {\bibinfo {volume} {87}},\ \bibinfo
  {pages} {86} (\bibinfo {year} {2015})}\BibitemShut {NoStop}%
\bibitem [{\citenamefont {Christodoulidi}\ \emph {et~al.}(2019)\citenamefont
  {Christodoulidi}, \citenamefont {Hone},\ and\ \citenamefont
  {Kouloukas}}]{C19}%
  \BibitemOpen
  \bibfield  {author} {\bibinfo {author} {\bibfnamefont {H.}~\bibnamefont
  {Christodoulidi}}, \bibinfo {author} {\bibfnamefont {A.~N.~W.}\ \bibnamefont
  {Hone}},\ and\ \bibinfo {author} {\bibfnamefont {T.~E.}\ \bibnamefont
  {Kouloukas}},\ }\href {https://doi.org/10.3934/jcd.2019011} {\bibinfo {title}
  {A new class of integrable lotka--volterra systems}} (\bibinfo {year}
  {2019})\BibitemShut {NoStop}%
\end{thebibliography}%

\paragraph{\textbf{Conserved quantities for graphs in \cref{fig:IntegrableGraphs}}}
\begin{description}
\item[a] \begin{align*}
f_1&=x_1+x_2+x_3\\
f_2&=x_1x_2x_3
\end{align*}

\item[b] \begin{align*}
f_1 &= x_1+x_2+x_3+x_4+x_5\\
f_3 &= x_1x_2x_3 + x_2x_3x_4 \\& \ \ \ \ \ \ \ \  + x_3x_4x_5 + x_4x_5x_1+x_5x_1x_2\\
f_5 &= x_1x_2x_3x_4x_5
\end{align*}

\item[c] \begin{align*}
f_1 &= x_1+x_2+x_3+x_4+x_5\\
f_3 &= (x_1+x_2+x_3)x_4x_5\\
f_9 &= x_1x_2x_3x_4^3x_5^3
\end{align*}

\item[d] \begin{align*}
f_1 &= x_1+x_2+x_3+x_4+x_5+x_6+x_7\\
f_3 &= x_1x_2x_3 + x_2x_3x_4 + x_3x_4x_5 + x_4x_5x_6 + x_5x_6x_7 \\ & \ \ \ \ x_6x_7x_1 + x_7x_1x_2 + x_1x_2x_5 + x_2x_3x_6 + x_3x_4x_7\\
& \ \ \ \ x_4x_5x_1+x_5x_6x_2+x_6x_7x_3 + x_7x_1x_4\\ 
f_5 &= x_1x_2x_3x_4x_5 + x_2x_3x_4x_5x_6 + x_3x_4x_5x_6x_7 \\ & \ \ \ \ +x_4x_5x_6x_7x_1 + x_5x_6x_7x_1x_2 + x_6x_7x_1x_2x_3  \\  \ \ \ \ \ \ \ \ & + x_7x_1x_2x_3x_4 \\
f_7 &= x_1x_2x_3x_4x_5x_6x_7
\end{align*}

\item[e] \begin{align*}
f_1 &= (x_1+x_2+x_3+x_4+x_5)+x_6+x_7\\
f_3 &= (x_1+x_2+x_3+x_4+x_5)x_6x_7\\
f_5 &= (x_1x_2x_3 + x_2x_3x_4  \\ \ \ \ \ & + x_3x_4x_5 + x_4x_5x_1+x_5x_1x_2)x_6^3x_7^3\\
f_{15} &= x_1x_2x_3x_4x_5x_6^5x_7^5
\end{align*}

\item[f] \begin{align*}
f_1 &= (x_1+x_2+x_3+x_4+x_5)+x_6+x_7\\
f_3 &= (x_1+x_2+x_3+x_4+x_5)x_6x_7\\
f_{9} &= ((x_1+x_2+x_3)x_4x_5)x_6^3x_7^3\\
f_{27} &= ((x_1x_2x_3)x_4^3x_5^3)x_6^9x_7^9
\end{align*}

\item[g] \begin{align*}
f_1 &= (x_1+x_2+x_3+x_4+x_5)+x_6+x_7\\
f_3 &= (x_1+x_2+x_3)x_4x_5 + x_4x_5x_6 \\ & \ \ \ \ \ \ \ \  + x_5x_6x_7 + x_6x_7(x_1+x_2+x_3) \\ & \ \ \ \ \ \ \ \  +x_7(x_1+x_2+x_3)x_4\\
f_{5} &= (x_1+x_2+x_3)x_4x_5x_6x_7\\
f_{15} &= x_1x_2x_3x_4^3x_5^3x_6^3x_7^3
\end{align*}

\item[h] \begin{align*}
f_1 &= (x_1+x_2+x_3)+(x_4+x_5+x_6)+x_7\\
f_3 &= (x_1+x_2+x_3)(x_4+x_5+x_6)x_7\\
f_{9} &= (x_1+x_2+x_3)^3(x_4x_5x_6)x_7^3 \\ & \ \ \ \ \ \ \ \  +(x_1x_2x_3)(x_4+x_5+x_6)^3x_7^3\\
f_{15} &= x_1x_2x_3x_4x_5x_6x_7^3
\end{align*}
\end{description}
\end{document}